\documentclass[preprint,12pt]{elsarticle}

\usepackage{graphics}
\usepackage{graphicx}
\usepackage{amssymb}
\usepackage{amsmath}
\usepackage{amsthm}

\usepackage{subfig}

\usepackage{booktabs}
\usepackage{multirow}
\usepackage{siunitx}
\usepackage{framed}
\usepackage[pagewise]{lineno}

\usepackage{ifthen,verbatim}
\usepackage{floatflt,graphicx,graphics,here}

\usepackage[urlcolor=blue,colorlinks=true]{hyperref}

\usepackage{tikz} 
\usepackage{mathrsfs}
\usepackage{amsmath,amssymb,amsthm}
\usepackage{mathtools}
\usepackage{pdfpages}
\usepackage{a4wide}
\usepackage{multirow}
\usepackage{url}
\usepackage{subfig}
\usepackage{booktabs}

\usepackage{interval}
\intervalconfig{soft open fences}

\numberwithin{equation}{section}
\nonstopmode
\numberwithin{equation}{section}
\theoremstyle{plain}

\theoremstyle{definition}

\theoremstyle{remark}

\newcommand{\R}{\mathbb{R}}

\newcommand{\B}{\mathbb{B}}

\newcommand{\uhp}{\mathbb{H}}

{\qed\bigskip}

\newcounter{alphabet}

\newcounter{minutes}
\divide\time by 60
\newcounter{hours}
\multiply\time by 60
\addtocounter{minutes}{-\time}

\usepackage{graphicx}
\graphicspath{{PDF/}{Graphics/}{Drafts/}}

\journal{}

\begin{document} 
\begin{frontmatter}



\title{Fast capacity computation for maze-like configurations}


\author[hh]{Harri Hakula}
\address[hh]{Aalto University\\
Department of Mathematics and System Analysis\\
P.O. Box 11100\\
FI--00076 Aalto, Finland
}
\ead{Harri.Hakula@aalto.fi}

\author[or]{Oona Rainio}
\address[or]{
University of Turku\\
Department of Mathematics and Statistics\\
FI--20014 University of Turku, Finland
}
\ead{ormrai@utu.fi}

\author[mv]{Matti Vuorinen}
\address[mv]{
University of Turku\\
Department of Mathematics and Statistics\\
FI--20014 University of Turku, Finland
}
\ead{vuorinen@utu.fi}


\begin{abstract}
We study the conformal capacity ${\rm cap}(\Omega,K)$ where $\Omega$ is a bounded domain of $\mathbb{R}^2$ and $K$ is a compact connected set in $\Omega$. Because the exact numerical value of the capacity is known only in a handful of special cases, it is important to find estimates for the capacity in terms of domain functionals, simpler than the capacity itself. Here, we study condensers of maze-like structure  and compute their  capacity by means of a high-order $hp$- finite element method. We compare these numerical results to the estimates given by the quasihyperbolic length and perimeter of the compact set. In particular, we consider the behaviour of these value pairs, numerical results and estimates, when the structure parameters vary and the walls of the maze approach the compact set. Over the configurations covered in the numerical 
experiments, the quasihyperbolic estimates are shown to have the desired asymptotic properties
and superior computational efficiencies once the case-specific analysis is completed.
\end{abstract}


\end{frontmatter}


\bibliographystyle{model1-num-names}

%


\footnotetext{\texttt{{\tiny File:~\jobname .tex, printed: \number\year-%
\number\month-\number\day, \thehours.\ifnum\theminutes<10{0}\fi\theminutes}}}
\makeatletter

\makeatother


\section{Introduction}
The design of many advanced engineering components, from microscopic sensors to 
high-performance electronic coolers, is invariably linked to a problem of geometry. 
The efficiency of devices like interdigitated capacitors or radial fin heat sinks depends directly on optimizing their complex, often ``maze-like'' physical structures. Of course, there exists vast number
of publications, and we can only cite some representative works \cite{hob,hus,kha}.
Although these applications seem distinct -- one concerning electrical fields and the other heat flow --
the underlying physical quantities of interest are related. From the point-of-view of applied
mathematics the core of the problem lies in solving Laplace's equation for a potential field, be it electric or thermal, within the domain defined by the component's boundaries. The solution to this is a harmonic potential function.

In this paper we study two-dimensional mazes of varying geometric complexity. The unifying
model is that of a \emph{condenser}. Formally, a condenser is a pair $(\Omega,K)$ consisting of a domain $\Omega\subset\R^n$ and its compact, non-empty subset $K\subset \Omega$. The walls of the maze are part of the boundary of the domain $\Omega$ and the compact set $K$ is a continuum which passes through the whole labyrinth like ``Ariadne's thread'' but carefully avoiding the walls.
The notion of condenser capacity is one of the key notions of mathematical physics related to the fields of electrostatistics and mathematical potential theory \cite{ps}, and has
many applications to engineering as noted above, see also \cite{sl}.
It is important to notice that the configuration we are focusing on 
seem to give nearly theoretically maximal capacity in terms of the quotient $m(\Omega)/d(K, \partial \Omega)^2$,
where $m(\Omega)$ is the measure and $d(\cdot,\cdot)$ a distance function.
The standard engineering approach is to employ some
numerical methods such as the finite element method (FEM). Here the reference results are
computed using a fine-tuned $hp$-FEM solver capable of handling curved boundaries and
optimal grading of the meshes in case of corner singularities. Fast, yet asymptotically
accurate approximations can be derived by employing techniques based on the quasihyperbolic
metric which can be considered natural in this setting. Similar ideas have been studied in
more constraint configurations before \cite{sa}. Another aspect is that in many 
applications the performance of the device does not require complex geometries, and indeed,
manufacturing may be constrained by the available tooling. We discuss a model problem
where the quantity of interest does not change significantly even though the underlying
configuration is simplified quite considerably. This, of course, can be interpreted as
an example of analysis-aware defeaturing, which is a highly topical subject \cite{bu,an}.

Condenser capacity is also important in geometric function theory \cite{du}.  
The capacity we study is the \emph{conformal capacity of a condenser}. 
The conformal capacity of a condenser  $(\Omega,K)$  is defined as \cite{hkv}
\begin{align}\label{def_condensercap}
{\rm cap}(\Omega,K)=\inf_{u\in A}\int_{\Omega}|\nabla u|^n \, dm,
\end{align}
where $A$ is the class of $C^{\infty}_0(\Omega)$ functions $u: \Omega \to \rinterval{0}{\infty} $
such that $u(x)\geq1$ for all $x \in K$ and $dm$ is the $n$-dimensional Lebesgue measure. 

Let us make three additional assumptions here: Let the dimension be $n=2$, suppose that the boundary components of $\Omega$ are non-degenerate continua, and choose the compact set $K$ so that it is connected. In this case, it is a basic fact that the infimum \eqref{def_condensercap} is in fact a minimum and the related extremal function $u$ is a harmonic
function, called as {\it the potential function} of the condenser.

The numerical value of the conformal capacity gives us information about the size of the set $K$ and its position with respect to the boundary of the domain $\Omega$. Finding this numerical value using analytical methods is, as a rule, impossible and therefore many methods have been developed to estimate
the values of condenser capacities such as isoperimetric inequalities \cite{ps}
and symmetrization methods \cite{bae}. Analytical methods are applicable to produce this numerical value only in special cases. In the case where $\Omega$ is the unit disk and $K$ is a connected set, there are
examples of conformal mappings of $\Omega \setminus K$ onto
canonical domains \cite{sl,ky}. 

The main novelty of this paper is that we not only demonstrate that with a fine-tuned $hp$-FEM solver
one can obtain highly reliable reference results even in very complex geometries but also
show that efficient approximations based on the quasihyperbolic length and perimeter are applicable.
These estimates are superior in terms of computational complexity. 
However, these advantages can only be realised after problem-specific calibration.

\subsection{Brief Review of Literature and Previous Work}
Soon after the introduction of the extremal length of a planar curve family and its reciprocal, the modulus of a curve family,
by Ahlfors and Beurling \cite{ab}, these notions became popular tools in the study of geometric function theory and 
quasiconformal mapping theory \cite{lv, du}. These popular tools were generalized by Fuglede \cite{f} to the case of $\mathbb{R}^n, n\ge 3,$
and this was particularly noteworthy because in higher dimensions much fewer effective tools were available than
in the planar case \cite{gmp,gor,hkm,hkv, vais} and in the
case of metric measure spaces \cite{hei}. It also turned out that the conformal capacity of a condenser
is related to the modulus: this capacity is equal to the modulus of the family of all curves joining $\partial \Omega$ with $K.$
The modulus of a curve family is thus connected with potential theory and
Ohtsuka \cite{o} continued Fuglede's work in this direction.

During the past half a century, parallel to the development of computer technology, there has been a lot of research on various methods to numerically compute conformal mappings and capacities \cite{avv,dali}. We refer the readers to the monographs \cite{dt,past,crow}, the comprehensive survey paper of R. Wegmann \cite{w}, and the short historical remarks in
\cite[pp. 14-16]{ps}, \cite[p. ix]{sl}, and \cite[pp. 8-12]{ky}.

G. P\'{o}lya  and G. Szeg\"{o} \cite{ps} developed systematically methods to estimate several domain functionals such as moment of inertia, principal frequency, torsion rigidity and, above all, condenser capacities. Their idea was to find estimates in terms of functionals simpler than the functional of interest. Thus, they used notions such as area and perimeter of sets to estimate capacities. They also
analyzed the sharpness of the results and noted that the extremal configuration often exhibited some symmetry. F.W. Gehring based his theory of quasiconformal mappings in space  \cite{gmp} on the conformal capacity, extending and applying the results
of  \cite{ps} to the case of conformal capacity. He also introduced domain functionals in terms of hyperbolic geometry which indeed is the natural setup in the case of conformal capacity \cite{g}. We, too, have applied numerical methods and hyperbolic geometry to study conformal capacity of condensers in a series of our earlier papers \cite{hrv,nrv1,nrv2,nv1,nv2}.

\subsection{Organization}
The structure of this paper is as follows. In Section 2, we introduce some basic notions from hyperbolic geometry and conformal mapping theory, and discuss the finite element method used in reference computations. In Section 3, we study two maze structures. The first of them is a subset of a rectangle where all the walls are segments parallel to the same side of the rectangle and the walls are evenly spaced. The second maze is a 
subset of the unit disk where the walls form a spiral like structure consisting of arcs of circles centered at the origin and radial segments joining these arcs. We study the behaviour of the capacity when the distance between the walls or the arcs decreases. In Section 4, our maze is a subset of an annulus and its walls are radial segments. In Section 5, we consider the case when $\Omega$ is the unit disk and the set $K$ is a union of disks of equal radii with centers at the same distance from origin and tangent to each other. It is natural to expect that a small part of the boundary of $K$, close to the points of contact of two disks is numerically negligible or "invisible" for the computation of the capacity of $(\Omega,K).$ We quantify this heuristic idea and present our results in the form of tables and graphics.

\section{Preliminaries}

In this section we first introduce some notations and conventions used in the sequel.
Next some basic notions from hyperbolic geometry and conformal mapping theory are covered, and finally, 
the finite element method ($hp$-FEM) used in reference computations is outlined 
including relevant details on error estimation and meshing.
\subsection{Notations and Conventions Employed}
For any point $x$ in a domain $\Omega\subsetneq\R^n$, we denote the Euclidean distance to the boundary $\partial \Omega$, defined as the infimum $\inf_{z\in\partial \Omega}|x-z|$, by $d_\Omega(x)$. We denote the $n$-dimensional $x$-centered open ball with radius $r>0$ by $B^n(x,r)=\{y\in\R^n\text{ : }|x-y|<r\}$ and the corresponding sphere by $S^{n-1}(x,r)=\{y\in\R^n\text{ : }|x-y|=r\}$. The argument of a non-zero complex number $z$ is denoted by ${\rm arg}(z)$ and the complex conjugate of $z$ by $\overline{z}$.



\subsection{Hyperbolic metric}
In the  upper half-space $\uhp^n=\{(x_1,...,x_n)\in\R^n\,:\,x_n>0\}$ and in the unit ball $\B^n=\{x\in\R^n\,:\,|x|<1\}$, the \emph{hyperbolic metric} can be defined by the formulas \cite[(4.8), p. 52 \& (4.14), p. 55]{hkv}, \cite{be}
\begin{align*}
\text{cosh}\rho_{\uhp^n}(x,y)&=1+\frac{|x-y|^2}{2d_{\uhp^n}(x)d_{\uhp^n}(y)},\quad x,y\in\uhp^n,\\
\text{sinh}^2\frac{\rho_{\B^n}(x,y)}{2}&=\frac{|x-y|^2}{(1-|x|^2)(1-|y|^2)},\quad x,y\in\B^n.
\end{align*}
In the special case $n=2$, we can simplify these formulas into
\begin{align*}
\text{tanh}\frac{\rho_{\uhp^2}(x,y)}{2}=\left|\frac{x-y}{x-\overline{y}}\right|,\quad
\text{tanh}\frac{\rho_{\B^2}(x,y)}{2}=\left|\frac{x-y}{1-x\overline{y}}\right|.
\end{align*}
\subsection{Quasihyperbolic metric} 
For a domain $\Omega\subsetneq\R^n$, the \emph{quasihyperbolic distance} between two points $x,y\in \Omega$ is \cite[(5.2), p. 68]{hkv} 
\begin{align*}
k_\Omega(x,y)=\inf_{\alpha\in\Gamma_{xy}}\int_\alpha w(u)|du|,   
\end{align*}
where $\Gamma_{xy}$ is the family of all rectifiable curves in $\Omega$ joining $x$ and $y$ and $w(u)$ is the weight function defined as \cite[(5.1), p. 68]{hkv}
\begin{align*}
w:\Omega\to (0,\infty),\quad w(u)=1\slash d_\Omega(u),\quad u\in \Omega.    
\end{align*}
If $\Omega$ is a simply-connected planar domain, then the quasihyperbolic metric fulfills the inequality
\begin{equation} \label{kvsrho}
\rho_{\Omega}(x,y)\leq2k_{\Omega}(x,y)\leq4\rho_{\Omega}(x,y)     
\end{equation}
for all points $x,y\in \Omega$ \cite[p.~21, (4.15)]{garmar}.



\subsection{Reference method for numerical solution of PDEs: $hp$-finite element method}\label{sec:hpfem} 
The $hp$-FEM is the method of choice for problems with strong corner singularities.
It can be shown that if the computational domain is properly discretized, exponential
convergence is achievable. This was first shown by Babu{\v{s}}ka and Guo \cite{BaGuo}.
For accessible treatment of the topic, we refer to Schwab \cite{schwab}.

Since the finite element methods are mature, different forms of error estimation are readily available.
One such method, the so-called auxiliary subspace method, has been used here to give us high confidence
in the numerical reference results (see \cite{hno}).

\begin{figure}
  \centering
  \subfloat{\includegraphics[width=0.45\textwidth]{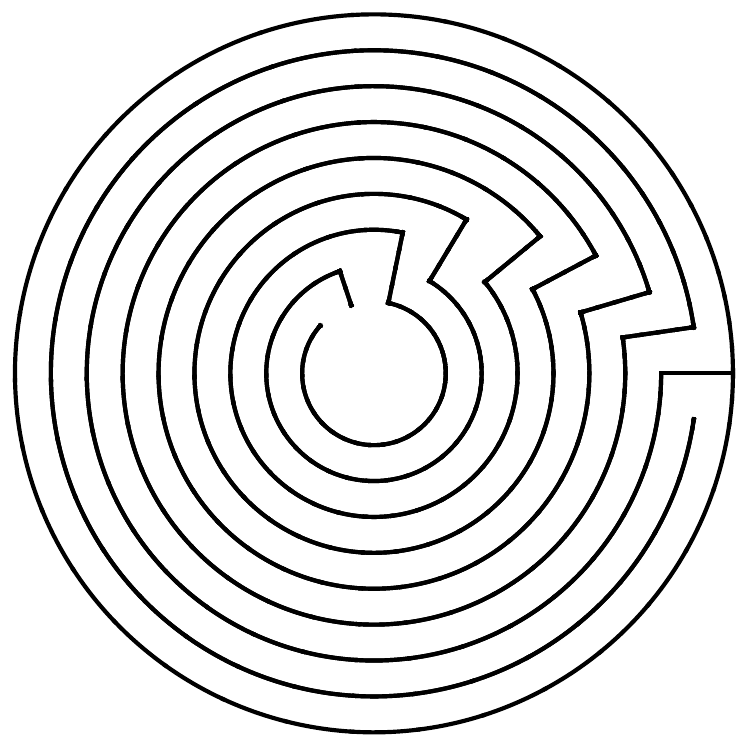}}\quad
  \subfloat{\includegraphics[width=0.45\textwidth]{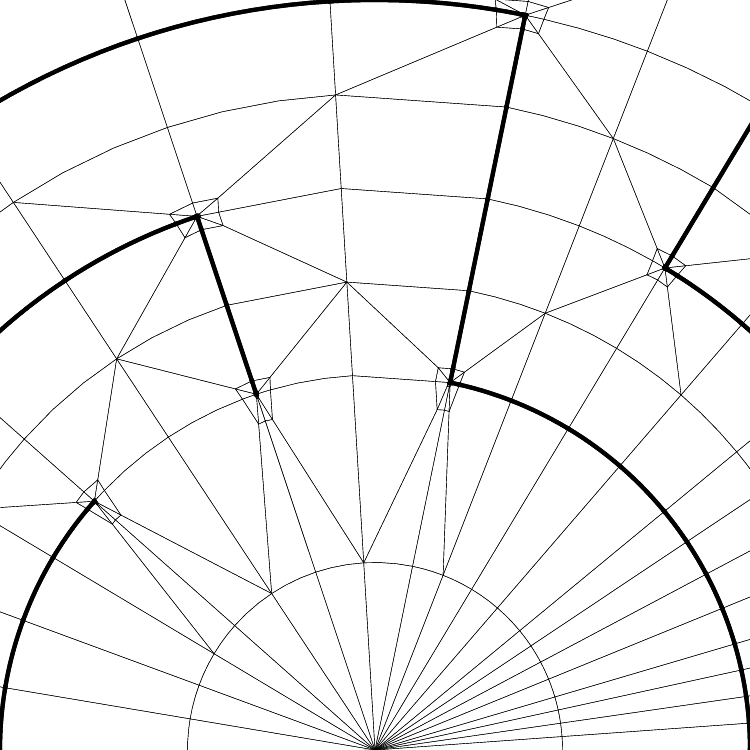}}
  \caption{Meshes for mazes. (Left) Circular maze. (Right) Detail
near the center, where the condenser boundaries are emphasized with thicker lines.
Strong mesh refinements are applied at every corner.
See~\cite{HaTu}.}\label{fig:refinement}
\end{figure}

An example of the meshing procedure is given in Figure~\ref{fig:refinement}.
The same approach has been applied to every corner of every maze considered below.

\section{Square and circular mazes}

Let us begin with maze-like configurations. 
First, the domain $\Omega$ is a square and the
walls are evenly spaced horizontal segments. 
Second, the domain is the unit disk and the walls are subarcs of circles 
centered at the origin with evenly spaces radii and these arcs are
joined by radial segments so that the walls form a spiral like connected set.
For short, we refer to these two mazes as the square maze and circular maze. The compact
set in both cases is the ``Ariadne's thread'' going through the whole maze and avoiding
the walls as much as possible.
 
\begin{figure}
  \centering
  \subfloat[{$m=11$}.]{
    \includegraphics[width=0.45\textwidth]{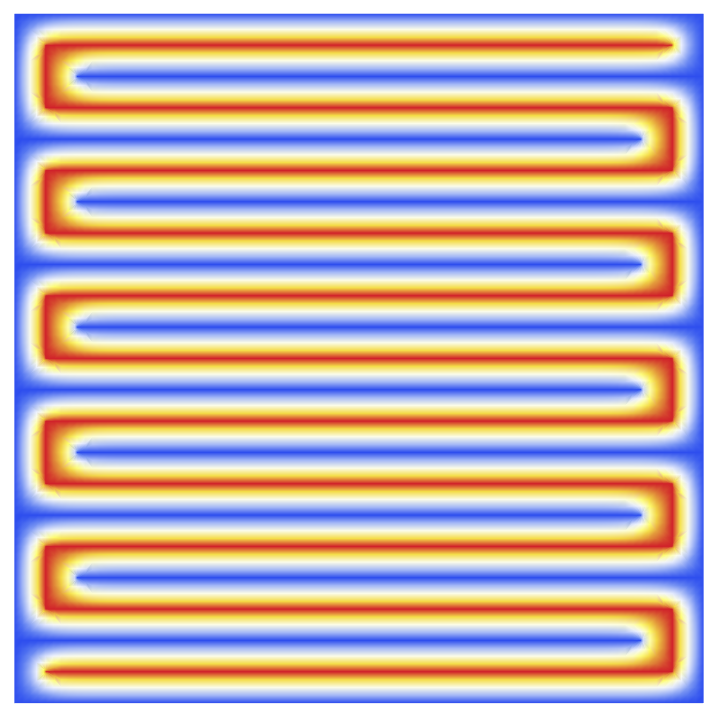}
  }
  \subfloat[{$m=14$}.]{
    \includegraphics[width=0.45\textwidth]{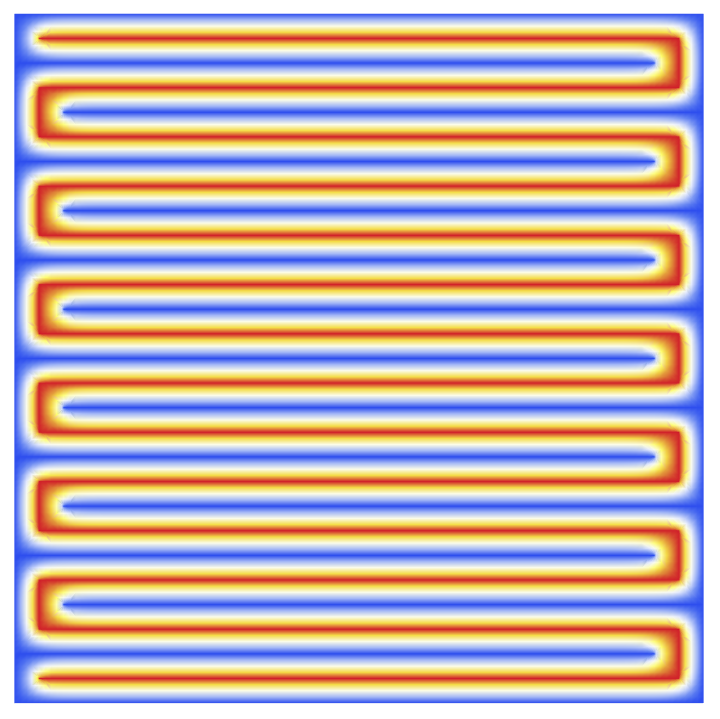}
  }
  \caption{Two potential functions of two square mazes.}
  \label{fig:comb}
\end{figure}

\subsection{Construction of the square maze} 
Our first type of condenser $(\Omega,K)$ is the square maze. For an integer $m\geq3$, the square maze domain $\Omega$ is a simply-connected subdomain of the unit square with additional $m-1$ horizontal Euclidean line segments inside to the unit square, referred here as spikes. These spikes start alternatively from the left and the right edges of the unit square so that the vertical distance between them is $1/m$. The Euclidean length of each spike is $1-1/m$. The first spike is $[i/m,1-1/m+i/m]$, the second spike $[1/m+2i/m,1+2i/m]$, the third spike $[3i/m,1-1/m+3i/m]$, and so on. The compact set $K$ is a polygonal chain circumventing the spikes inside the unit square. Examples can be seen in Figure \ref{fig:comb}.

As noted in \cite{nrv2}, the hyperbolic perimeter of the compact set can be used to estimate the condenser capacity. Since finding the exact value of the hyperbolic metric itself is complicated in our maze domain, we use the quasihyperbolic metric instead.

\subsection{Quasihyperbolic perimeter of the compact set in the square maze} 
  The compact set is a polygonal chain with $m$ horizontal line segments and $m-1$ vertical line segments inside the maze. The Euclidean length of the horizontal line segment is $1-1/m$ and the Euclidean distance from a point of this line segment to the domain boundary is $1/(2m)$ for all points of the given horizontal segment. The quasihyperbolic length of a single horizontal segment is therefore
\begin{align*}
\int^{1-1/m}_{0}\frac{1}{1/(2m)}du
=\int^{1-1/m}_{0}2mdu
=2m(1-1/m)
=2(m-1).
\end{align*}
The Euclidean length of the vertical line segment is $1/m$ and the Euclidean distance from its any point to the boundary $1/(2m)$, so the quasihyperbolic length of one vertical segment is
\begin{align*}
\int^{1/m}_{0}\frac{1}{1/(2m)}du=\int^{1/m}_{0}2mdu=2.
\end{align*}
Consequently, the quasihyperbolic length of the whole polygonal chain with $m$ horizontal and $m-1$ vertical segments is 
\begin{align*}
m\cdot2(m-1)+(m-1)\cdot2=2(m^2-1).
\end{align*}
The quasihyperbolic perimeter of the polygonal chain is twice of its quasihyperbolic length, giving us the final expression
\begin{align*}
4(m^2-1).    
\end{align*}

\begin{figure}
  \centering
  \subfloat[Error estimates.]{
    \includegraphics[width=0.45\textwidth]{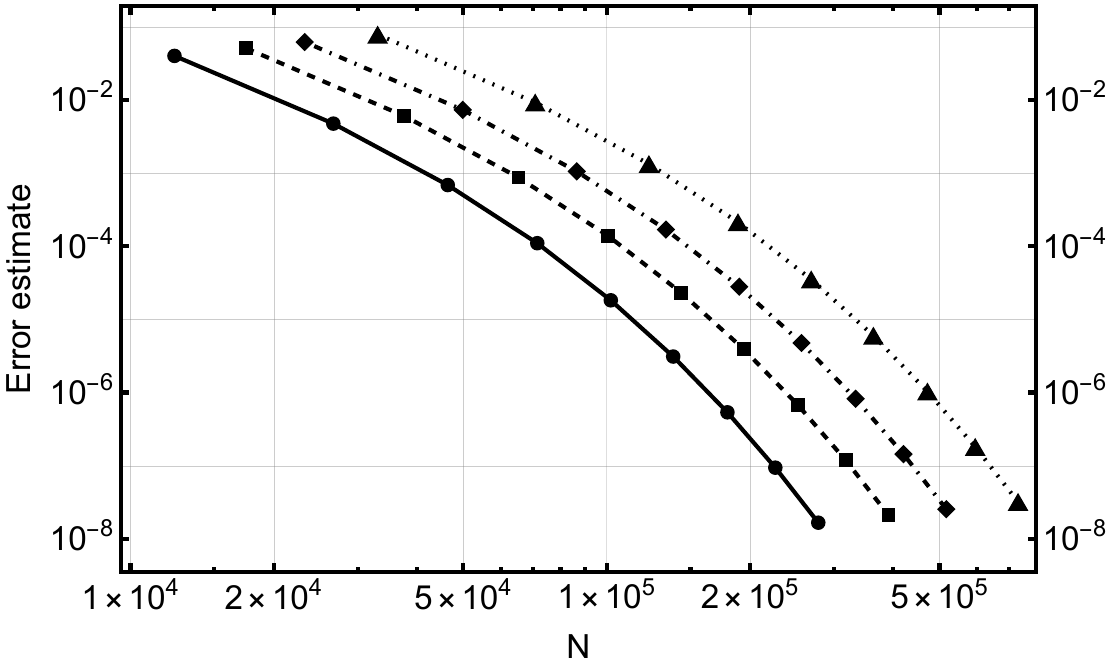}
  }
  \subfloat[Parameter dependence.]{
    \includegraphics[width=0.45\textwidth]{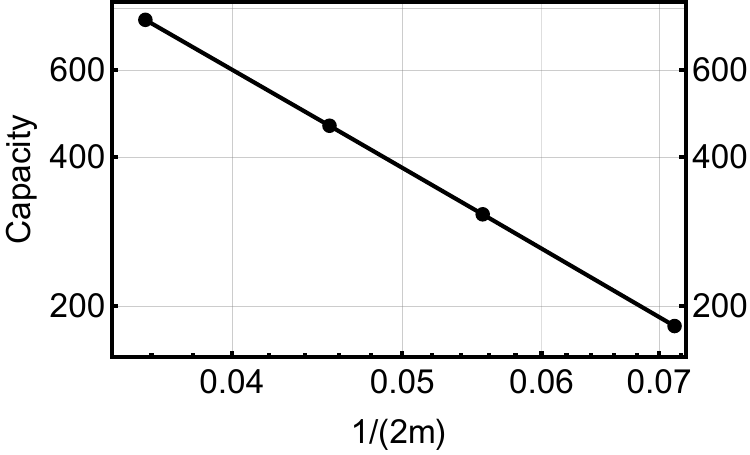}
  }
  \caption{Square maze: Error estimates and dependence of capacity on the parameter $m$
  which indicates the number of horizonal segments (hence the distance
  between segments is $1/m$) (loglog-plots). 
  (a) Error estimates of the capacity are shown as functions of the number of degrees of freedom $N$ in
 the $hp$-FEM solution. The four graphs represent parameter values $m = 7,9,11,14$ appearing from left to right in the figure. The points on the graphs correspond to $p=2,\ldots,10$ on a given mesh.
  (b) Dependence of capacity on the parameter $1/(2m)$ has the rate $\sim 2.05$.
  }
  \label{fig:combconvergence}
\end{figure}

\begin{table}
    \centering
    \caption{Quasihyperbolic length and quasihyperbolic perimeter of the square maze for different numbers $m$ of arcs in the spiral domain. 
    Computed reference capacities are at $p=10$.
    The error estimate $\varepsilon$ is given in compressed form. $N$ is the number of degrees of freedom.}
    \label{tbl:sqmaze}
    \begin{tabular}{rrrlrr}
        \toprule
\multirow{2}[3]{*}{$m$} & \multicolumn{2}{c}{Quasihyperbolic} & \multicolumn{3}{c}{$hp$-FEM $p=10$} \\
        \cmidrule(lr){2-3} \cmidrule(lr){4-6}
                   & Length & Perimeter &    Capacity & $|\lceil\log_{10}(\varepsilon)\rceil|$ & $N$               \\
        \midrule
        7  &  96 & 192 & 182.21381864123760 & 7 & 278160 \\
        9  & 160 & 320 & 306.32589930610663 & 7 & 390720 \\
        11 & 240 & 480 & 462.43797997097550 & 7 & 516720 \\
        14 & 390 & 780 & 756.60610096827860 & 7 & 730920 \\
        \bottomrule
    \end{tabular}
\end{table}

\subsection{Construction of the circular maze} 
\begin{figure}
  \centering
  \subfloat[{$M=5$.}]{
    \includegraphics[width=0.45\textwidth]{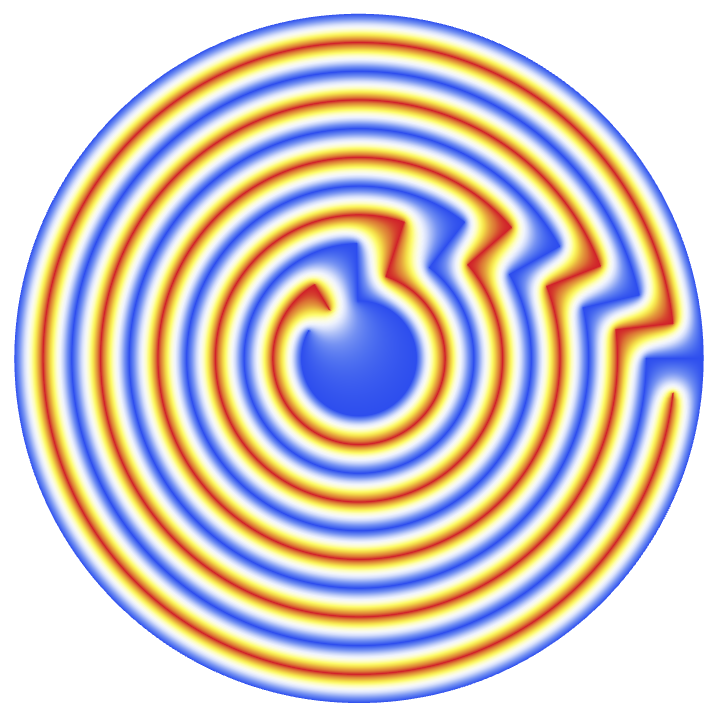}
  }
  \subfloat[{$M=15$.}]{
    \includegraphics[width=0.45\textwidth]{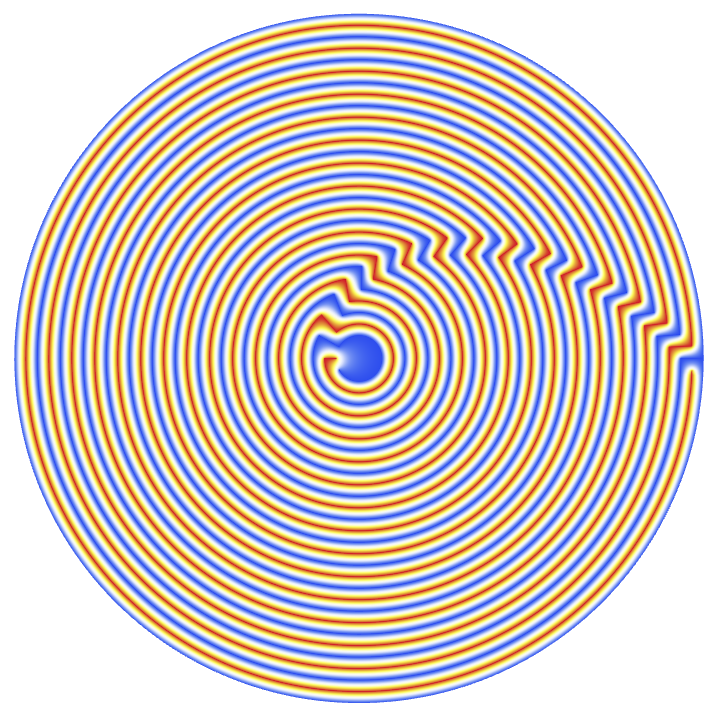}
  }
  \caption{Circular Maze: Two realisations.} 
  \label{fig:circularMazePotential}
\end{figure}

A circular maze condenser $(\Omega,K)$ consists of a simply-connected subdomain $\Omega$ of the unit disk $\B^2$ and a compact set $K$. The boundary of $\Omega$ has two parts, the unit circle and a spiral-like curve that forms the walls of the maze. The set $K$, ``Ariadne's thread'', travels through the whole maze staying as far from the walls as possible. For a given integer $m\geq3$, the walls are formed by sub-arcs of $S^1(0,r_k)$, $r_k=(m-k)/m$, $k=1,2,...,m-1$, and radial segments connecting sub-arcs of $S^1(0,r_k)$ and $S^1(0,r_{k+1})$. Similarly, $K$ is formed by sub-arcs of $S^1(0,r_k+1/(2m))$ and connecting radial segments.
For illustrations, see Figure~\ref{fig:circularMazePotential}.

Let us describe the structure of the walls and the set $K$ in detail. The walls consist of the sub-arcs
\begin{align*}
\{z\,:\,|z|=r_k,\,{\rm arg}(z)\notin(\theta_k,\theta_k+\alpha_k)\}    
\end{align*}
of $S^1(0,r_k)$, where 
\begin{align*}
\alpha_k=2{\rm arcsin}\left(\frac{1}{2(m-k)-1}\right),\,k=1,2,...,m-2;\quad
\alpha_{m-1}=2{\rm arcsin}(1/2),
\end{align*}
$\theta_{k+1}=\theta_k+\alpha_{k+1}$, and the radial segment is
\begin{align*}
\{z\,:\,{\rm arg}(z)=\theta_k,\,r_{k+1}\leq|z|\leq r_k\}.     
\end{align*}
The set $K$ consists of the sub-arcs
\begin{align*}
\left\{\left(r_k+\frac{1}{2m}\right)e^{i\varphi}\,:\,\varphi\notin\left[\theta_k-\frac{\alpha_{k-1}}{2},\theta_k+\frac{\alpha_{k-1}}{2}\right]\right\}    
\end{align*}
and the radial segments
\begin{align*}
\left\{re^{i(\theta_k+\alpha_k/2)}\,:\,r_k-\frac{1}{2m}\leq r\leq r_k+\frac{1}{2m}\right\},   
\end{align*}
though the innermost radial segment is an exception to this. Thus, the points of $K$ are approximately at distance $1/(2m)$ to the boundary set $\partial \Omega$ of $\Omega$.

\subsection{Quasihyperbolic length of the compact set in the circular maze}
The compact set contains $m-1$ circular arcs and $m-1$ radial line segments.

The radius of the $k$th outermost arc, $k=1,...,m-1$, is $(m-k+1/2)/m$ and the central angle is $2\pi-\alpha_k$. For all $k=1,...,m-1$, the Euclidean length of the circle arcs is
\begin{align*}
(2\pi-\alpha_k)\left(r_k+\frac{1}{2m}\right)
=(2\pi-\alpha_k)(m-k+1/2)/m.
\end{align*}
The minimum distance from their any point to the the domain boundary is always $1/(2m)$. Consequently, the quasihyperbolic length of a single arc is
\begin{align*}
&\int^{(2\pi-\alpha_k)(m-k+1/2)/m}_02m\,du=(2\pi-\alpha_k)(2m-2k+1).
\end{align*}
The sum of the lengths of all the $m-1$ arcs is equal to
\begin{align*}
&\sum^{m-1}_{k=1}(2\pi-\alpha_k)(2m-2k+1).
\end{align*}

Each of the $m-2$ outermost radial line segments has a Euclidean length equal to $1/m$ and the innermost has a Euclidean length equal to $1/(2m)$. For such a point on a radial segment whose absolute value is $u$ with $r_k-1/(2m)\leq u\leq r_k+1/(2m)$, the distance to boundary is $u\sin(\alpha_k/2)$. The quasihyperbolic length of the $k$th outermost radial line segment, $k=1,...,m-2$, equals
\begin{align*}
&\int^{r_k+1/(2m)}_{r_k-1/(2m)}\frac{1}{u\sin(\alpha_k/2)}\,du  =\frac{1}{\sin(\alpha_k/2)}\int^{r_k+1/(2m)}_{r_k-1/(2m)}\frac{1}{u}\,du\\
&=\frac{1}{\sin(\alpha_k/2)}(\log|r_k+1/(2m)|-\log|r_k-1/(2m)|)\\
&=\frac{1}{\sin(\alpha_k/2)}(\log|m-k+1/2|-\log|m-k-1/2|)
\end{align*}
The innermost radial segment with $k=m-1$ has a quasihyperbolic length equal to
\begin{align*}
&\int^{r_{m-1}+1/(2m)}_{r_{m-1}}\frac{1}{u\sin(\alpha_k/2)}\,du 
=\frac{1}{\sin(\alpha_{m-1}/2)}(\log|r_{m-1}+1/(2m)|-\log|r_{m-1}|)\\
&=\frac{1}{\sin(\alpha_{m-1}/2)}(\log|3/2|)
\end{align*}
The sum of the quasihyperbolic lengths of the all the radial segments is equal to
\begin{align*}
\frac{\log|3/2|}{\sin(\alpha_{m-1}/2)}+\sum^{m-2}_{k=1}\left(\frac{1}{\sin(\alpha_k/2)}(\log|m-k+1/2|-\log|m-k-1/2|)\right).   
\end{align*}

Now, the quasihyperbolic length of the whole spiral is 
\begin{align*}
&\sum^{m-1}_{k=1}(2\pi-\alpha_k)(2m-2k+1)\\
&+\frac{\log|3/2|}{\sin(\alpha_{m-1}/2)}+\sum^{m-2}_{k=1}\left(\frac{1}{\sin(\alpha_k/2)}(\log|m-k+1/2|-\log|m-k-1/2|)\right)
\end{align*}
and its quasihyperbolic perimeter is twice of this. See Table~\ref{t4} for the quasihyperbolic lengths and perimeters for a few example values of $m$.
 
\begin{table}
    \centering
    \caption{Quasihyperbolic length and quasihyperbolic perimeter of the circular maze for different numbers $m$ of arcs in the spiral domain. 
    Computed reference capacities are at $p=10$.
    The error estimate $\varepsilon$ is given in compressed form. $N$ is the number of degrees of freedom.}
    \label{t4}
    \begin{tabular}{rrrlrr}
            \toprule
    \multirow{2}[3]{*}{$m$} & \multicolumn{2}{c}{Quasihyperbolic} & \multicolumn{3}{c}{$hp$-FEM $p=10$} \\
            \cmidrule(lr){2-3} \cmidrule(lr){4-6}
                       & Length & Perimeter &    Capacity & $|\lceil\log_{10}(\varepsilon)\rceil|$ & $N$               \\
            \midrule
        5  &  144 &  289  & \phantom{0}287.0934278396659 & 7 & 241670 \\
        7  &  294 &  588  & \phantom{0}583.6685510819680 & 7 & 373710 \\
        10 &  613 & 1227  & 1218.3273371997180 & 7 & 596980 \\
        12 &  889 & 1778  & 1766.2460049187985 & 6 & 762610 \\
        15 & 1397 & 2795  & 2777.8399104479818 & 5 & 1036270 \\
        \bottomrule
    \end{tabular}
\end{table}

\begin{figure}
  \centering
  \subfloat[Error estimates.]{
    \includegraphics[width=0.45\textwidth]{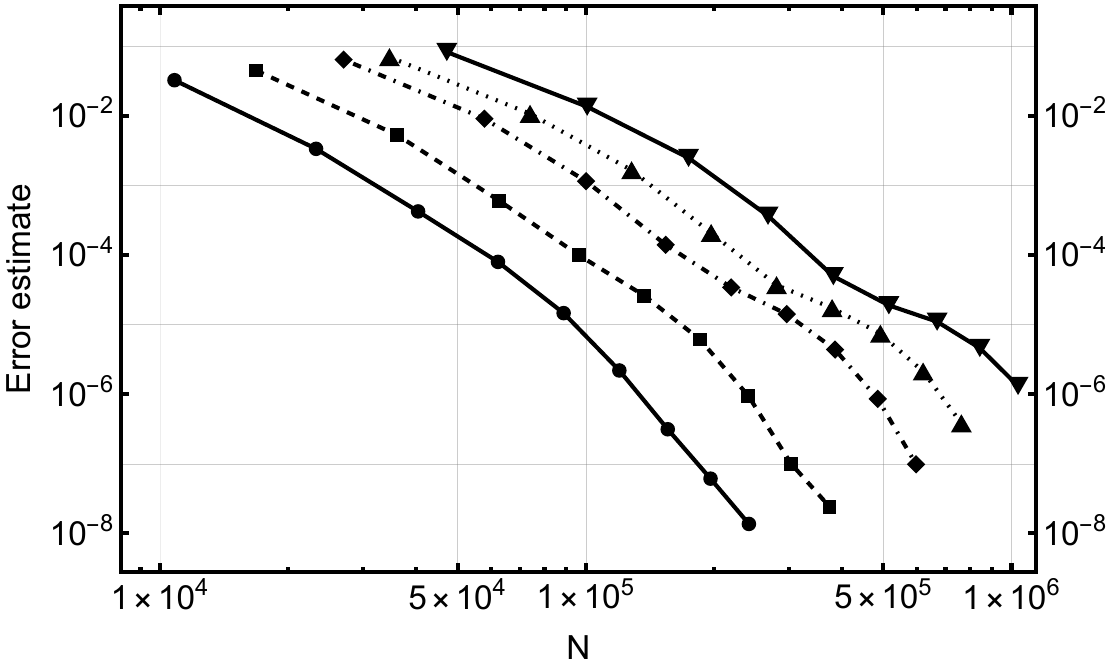}
  }
  \subfloat[Parameter dependence.]{
    \includegraphics[width=0.45\textwidth]{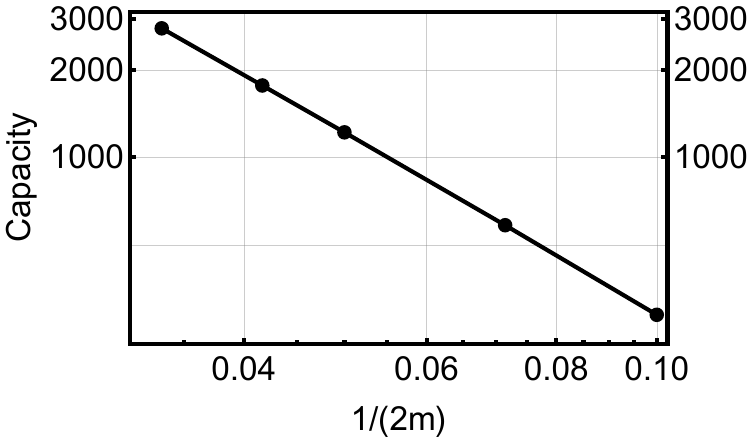}
  }
  \caption{Circular Maze: Error estimates and dependence of capacity on the parameter $m$  (loglog-plots). 
   (a) Error estimates in capacity are shown as functions of the number of degrees of freedom $N$ in
   the $hp$-FEM solution. The four graphs represent parameter values $m = 5,7,10,12,15$ appearing from left to right in the figure. The points on the graphs correspond to $p=2,\ldots,10$ on a given mesh. (b) Dependence of capacity on the parameter $1/(2m)$ has the rate $ \sim 2.06$.}
  \label{fig:circularconvergence}
\end{figure}

\begin{figure}
  \centering
    \includegraphics[width=0.45\textwidth]{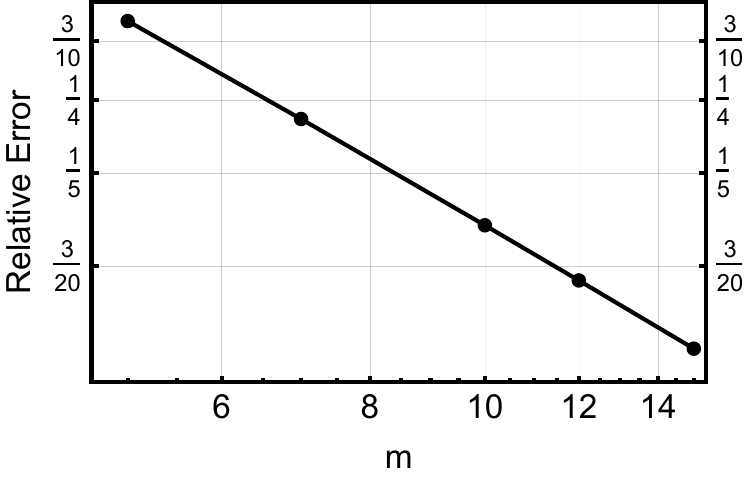}
  \caption{Circular Maze: Convergence of the relative error of the quasihyperbolic perimeter. The rate is $\sim 0.92$. (loglog-plot)}\label{fig:circularMazeConv}
\end{figure}

\subsection{Analysis of results: Conformance with theoretical predictions}
The results are collected in Tables~\ref{tbl:sqmaze} and \ref{t4}. The reference results are
illustrated in Figures~\ref{fig:combconvergence} and \ref{fig:circularconvergence}.
The relative error of the the quasihyperbolic perimeter as a function of the parameter $m$
is shown in Figure~\ref{fig:circularMazeConv}.

The $hp$-FEM solutions converge exponentially and the estimated errors indicate high level
of reliability. Moreover, the dependence of capacity on the parameter $1/(2m)$ has the rate $ \sim 2$
in both cases in conformance with the theoretical prediction \cite[Thm 7.1, pg. 20]{vais}.
It is evident that the quasihyperbolic perimeter is the natural metric in this setting.
As the number of segments increase, the estimate becomes asymptotically exact.
\section{Annulus with spikes}
Our next example considers the maze inside an annulus. As in the previous example, the
quasihyperbolic length and perimeter are computed in closed form and compared with the
$hp$-FEM results.
\subsection{Construction of the spiked annulus maze}
\begin{figure}
  \centering
  \subfloat[$M=10$]{
    \includegraphics[width=0.3\textwidth]{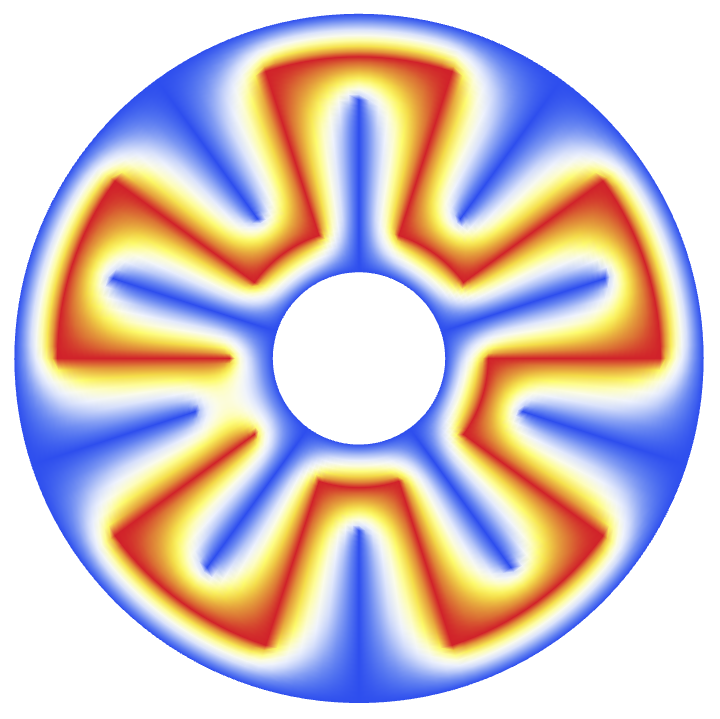}
  }
  \subfloat[$M=16$]{
    \includegraphics[width=0.3\textwidth]{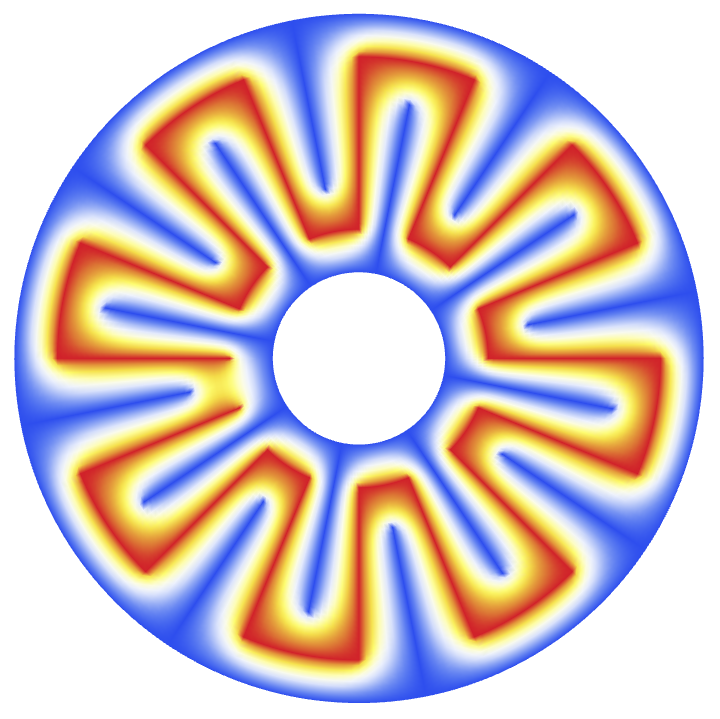}
  }
  \subfloat[$M=20$]{
    \includegraphics[width=0.3\textwidth]{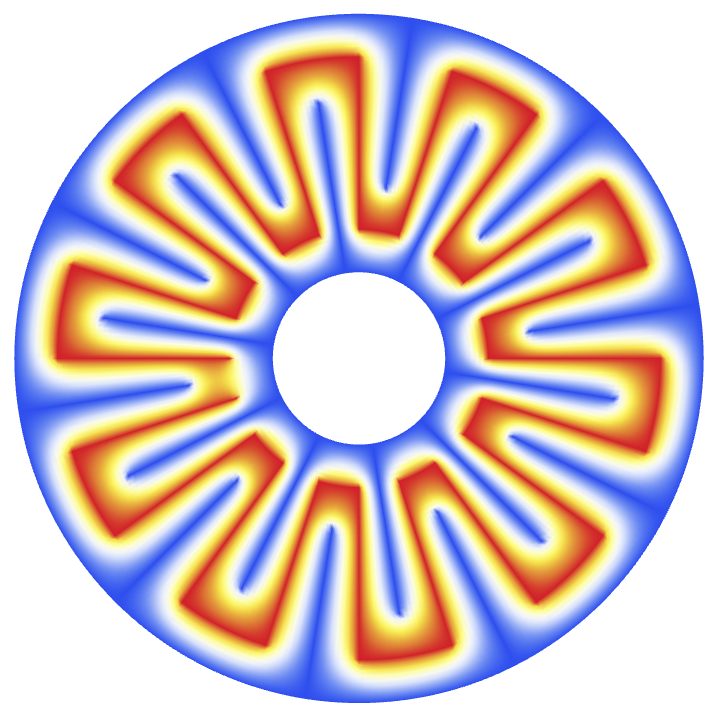}
  }
  \caption{Spiked Annulus mazes: The number of spikes is $M$, the radius of the inner boundary $r=1/4$, the length of the spikes $l=1/2$.}\label{fig:spikedannulusmaze}
\end{figure}

\begin{figure}
  \centering
  \subfloat{
    \includegraphics[width=0.45\textwidth]{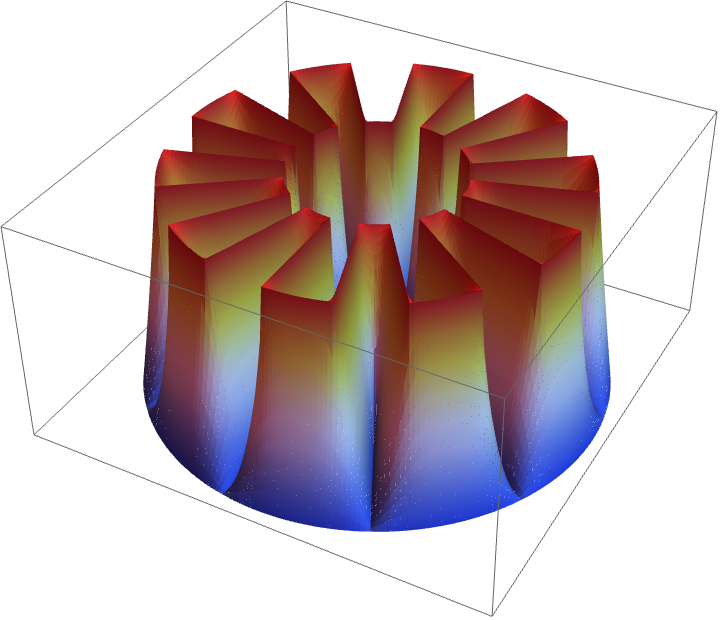}
  }
  \subfloat{
    \includegraphics[width=0.45\textwidth]{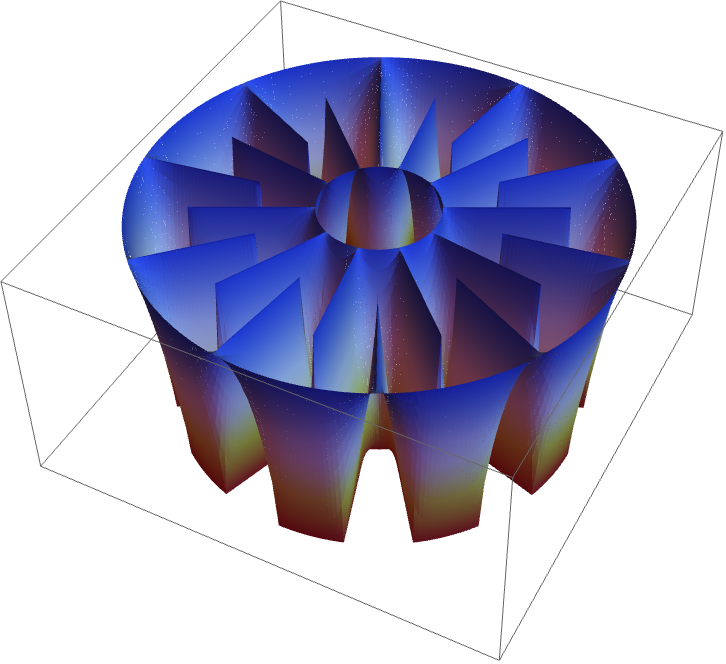}
  }
  \caption{Spiked Annulus mazes ($M=16$): Two perspectives of the potential function.}\label{fig:spikedannulusmaze3d}
\end{figure}

Annulus mazes are parameterised as given in Table~\ref{tbl:mazeparameters}.
Realisations are computed with $r_0=1/4$, $l_0=l_1=1/2$, and $r_1=1$,
and hence $R_0 = 3/8$ and $R_1=7/8$ resulting in $C_1 = 1/2$.
Every circular arc covers the sector of $2\pi/M$ radians, and
by construction $k_0 = (M/2)-1$ and $k_1 = M/2$.
Some realizations are illustrated in Figures~\ref{fig:spikedannulusmaze}
and \ref{fig:spikedannulusmaze3d}.

\begin{table}
\centering
\caption{Spiked annulus maze: Parameters.}\label{tbl:mazeparameters}.
\begin{tabular}{ll}\toprule
Parameter & Value \\ \midrule
$M$ & Total number of spikes (even) in the domain\\
$r_0$ & Inner radius of the domain\\
$r_1$ & Outer radius of the domain\\
$l_0$ & Length of inner spikes in the domain\\
$l_1$ & Length of outer spikes in the domain\\
$k_0$ & Number of inner arcs in the compact set\\
$k_1$ & Number of outer arcs in the compact set\\ 
$R_0$ & Inner arc radius in the compact set\\
$R_1$ & Outer arc radius in the compact set\\
$C_1$ & Length of the radial segment in the compact set\\
\bottomrule
\end{tabular}
\end{table}

\subsection{Quasihyperbolic length of the chain in the spiked annulus domain}
Suppose the annulus has $M$ spikes where $M\geq6,8,10,...$.
The chain consists of $M$ radial line segments, $(M/2-1)$ arcs with radius of $3/8$, and $M/2$ arcs with radius of $7/8$.
The Euclidean length of each radial segment is $1/2$. For $0\leq u\leq1/2$, the distance to the inner boundary of the domain is $u+1/8$, the distance to the outer boundary of the domain is $1/2-u+1/8=5/8-u$, and the distance the closest spike is $\sin(\pi/(2M))(u+3/8)$. For $M\geq6$, $\sin(\pi/(2M))<0.26$, the spike is the closest part of the boundary from the line segment. The quasihyperbolic length of the line segment is
\begin{align*}
&\int^{1/2}_0\frac{1}{\sin(\pi/(2M))(u+3/8)}du
=\frac{1}{\sin(\pi/(2M))}\int^{1/2}_0\frac{1}{u+3/8}du\\
&=\frac{1}{\sin(\pi/(2M))}(\log|8(1/2)+3|-\log|8\cdot0+3|)
=\frac{\log(7)-\log(3)}{\sin(\pi/(2M))}.
\end{align*}

A single piece of the arc with radius of $3/8$ has the Euclidean length of $3/8\cdot2\pi/M=3\pi/(4M)$ and a distance of $1/8$ to the boundary, so its quasihyperbolic length is
\begin{align*}
&\int^{3\pi/(4M)}_0\frac{1}{1/8}du
=\int^{3\pi/(4M)}_08du
=8\cdot3\pi/(4M)-8\cdot0
=6\pi/M.
\end{align*}

Similarly, a single piece of the arc with radius of $7/8$ has the Euclidean length of $7/8\cdot2\pi/M=7\pi/(4M)$ and a distance of $1/8$ to the boundary, so its quasihyperbolic length is
\begin{align*}
&\int^{7\pi/(4M)}_0\frac{1}{1/8}du
=\int^{7\pi/(4M)}_08du
=8\cdot7\pi/(4M)-8\cdot0
=14\pi/M.
\end{align*}

The quasihyperbolic length of the whole chain is
\begin{align*}
&\frac{M(\log(7)-\log(3))}{\sin(\pi/(2M))}+(M/2-1)\cdot 6\pi/M+M/2\cdot 14\pi/M\\
&=\frac{M(\log(7)-\log(3))}{\sin(\pi/(2M))}+3\pi(M-2)/M+7\pi
\end{align*}
and the quasihyperbolic perimeter is twice of this. See Table~\ref{t_qhp}.

\begin{table}
    \centering
    \caption{Quasihyperbolic length and quasihyperbolic perimeter of the chain in the spiked annulus domain for different numbers of spikes $M$. 
    Computed reference capacities are at $p=12$.
    The error estimate $\varepsilon$ is given in compressed form. $N$ is the number of degrees of freedom.}
    \label{t_qhp}
    \begin{tabular}{rrrlrr}
            \toprule
    \multirow{2}[3]{*}{$M$} & \multicolumn{2}{c}{Quasihyperbolic} & \multicolumn{3}{c}{$hp$-FEM $p=12$} \\
            \cmidrule(lr){2-3} \cmidrule(lr){4-6}
                       & Length & Perimeter &    Capacity & $|\lceil\log_{10}(\varepsilon)\rceil|$ & $N$               \\
            \midrule
        10 & 83  & 167 & \phantom{0}93.43853840480526   & 9 & 493583 \\
        16 & 168 & 337 & 167.62306400419487  & 8 & 789935 \\
        20 & 246 & 492 & 233.82949214101268  & 8 & 987503 \\
        \bottomrule
    \end{tabular}
\end{table}

\begin{figure}
  \centering
  \subfloat[Error estimates.]{
    \includegraphics[width=0.45\textwidth]{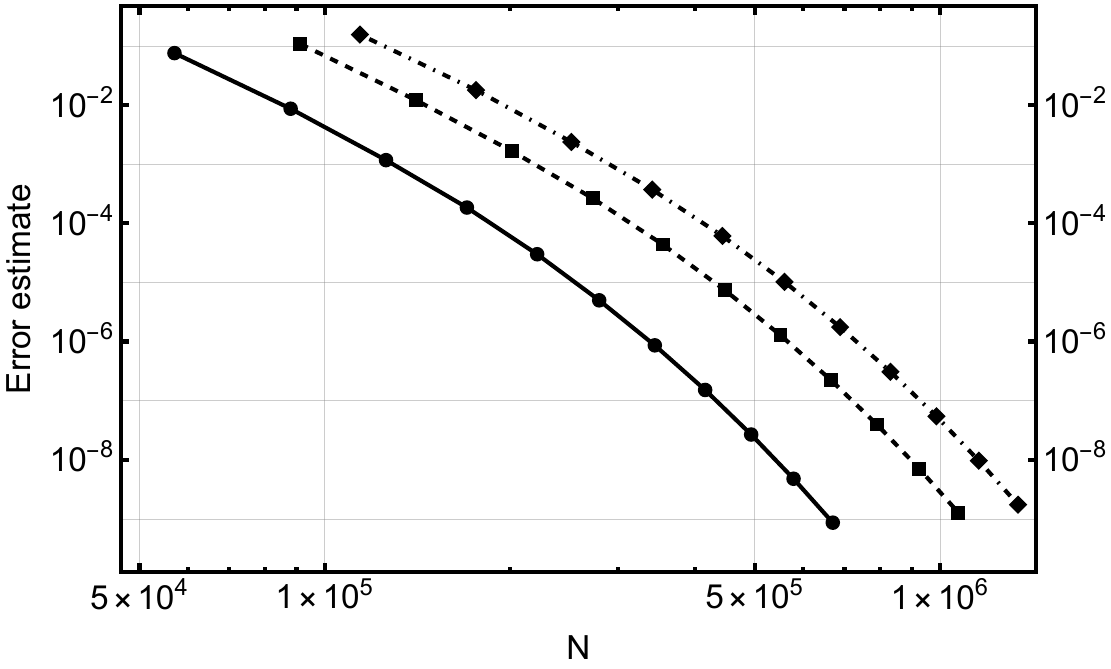}
  }
  \subfloat[Parameter dependence.]{\label{fig:spikedconvergenceb}
    \includegraphics[width=0.45\textwidth]{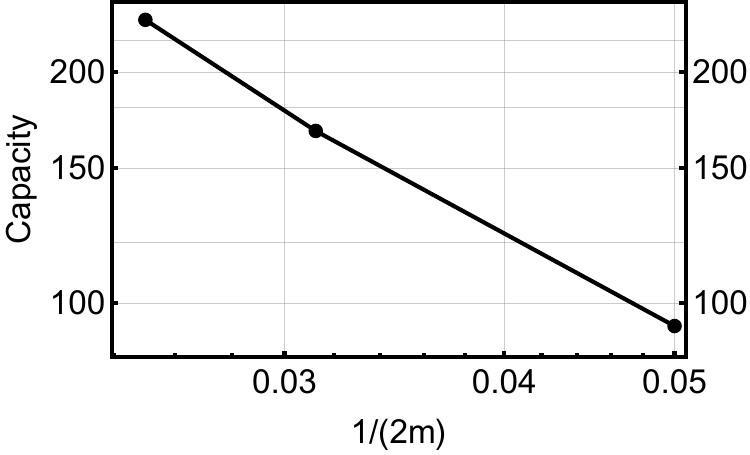}
  }
  \caption{Spiked Annulus: Error estimates and dependence of capacity on the parameter $m$  (loglog-plots). 
   (a) Error estimates in capacity are shown as functions of the number of degrees of freedom $N$ in
   the $hp$-FEM solution. The four graphs represent parameter values $m = 10,16,20$ appearing from left to right in the figure. The points on the graphs correspond to $p=2,\ldots,12$ on a given mesh. (b) Dependence of capacity on the parameter $1/(2M)$ has the rate $\sim 1.31$.}
  \label{fig:spikedconvergence}
\end{figure}

\subsection{Analysis of results}
The results are presented in Table~\ref{t_qhp} and Figure~\ref{fig:spikedconvergence}.
Again, the $hp$-FEM results are in line with expectations. 
Now the quasihyperbolic length is the natural choice for estimation.
However, the agreement is not as good as in the mazes above. Having said that, there
are no predictions for the theoretical dependence on the parameter $M$.
For the $hp$-FEM we observe a rate $\sim 1.3 < 2$. It appears that new theory is required to
verify the results.

\subsection{Comment on computational complexity}
The $hp$-FEM results have been computed on modern Apple Silicon Macs with 128GB RAM.
The assembly of the matrices, which is often the most expensive part in high-order FEM, 
has been parallelized. Only direct solvers provided by Mathematica have been used.

In this class of problems with complex geometries, the real cost is in defining the
geometry and the corresponding mesh. The high reliability obtained here is due to exact
geometry representation. In the absence of standard meshing tools supporting high-order FEM,
this requires manual intervention at every level, and makes it very difficult to
give a clear picture of the true costs involved.

In Table~\ref{tbl:timing} timing data on spiked annulus cases is given. Here the assembly has been computed
sequentially, but linear solver is the built-in one. As one would expect, the time for assembly,
or in other words, numerical integration in this context, grows linearly in $M$, but the solution times
grow slightly faster, yet are only a fraction of the time spent on assembly.
Also, as has been observed before,
the relative cost of auxiliary space error estimation decreases at the polynomial order increases.

\begin{table}
    \centering
    \caption{Spiked annulus domains: Timing data for $hp$-FEM. Breakdown for different parts of the solution process for $M$ and polynomial order $p$. $N_{\text{tot}}$ is the dimension of the fully integrated system including the auxiliary space degrees of freedom.}
    \label{tbl:timing}
    \begin{tabular}{rrrrrr}
            \toprule
        $M$ & $p$ & Assembly & Solve & Estimation & $N_{\text{tot}}$ \\
            \midrule
        10 & 8  & 36  & 2.0  & 1.6 & 344119 \\
           & 12 & 77  & 7.8  & 4.4 & 669927 \\
        16 & 8  & 54  & 3.2  & 2.7 & 550759 \\
           & 12 & 121 & 12.1 & 6.9 & 1072119 \\
        20 & 8  & 67  & 4.1  & 3.0 & 688519 \\
           & 12 & 155 & 15.4 & 8.9 & 1340247 \\
        \bottomrule
    \end{tabular}
\end{table}


\section{Non-visible parts of a compact set}

\begin{figure}
  \centering
  \subfloat[]{\label{fig:nonvisibleA}
    \includegraphics[width=0.45\textwidth]{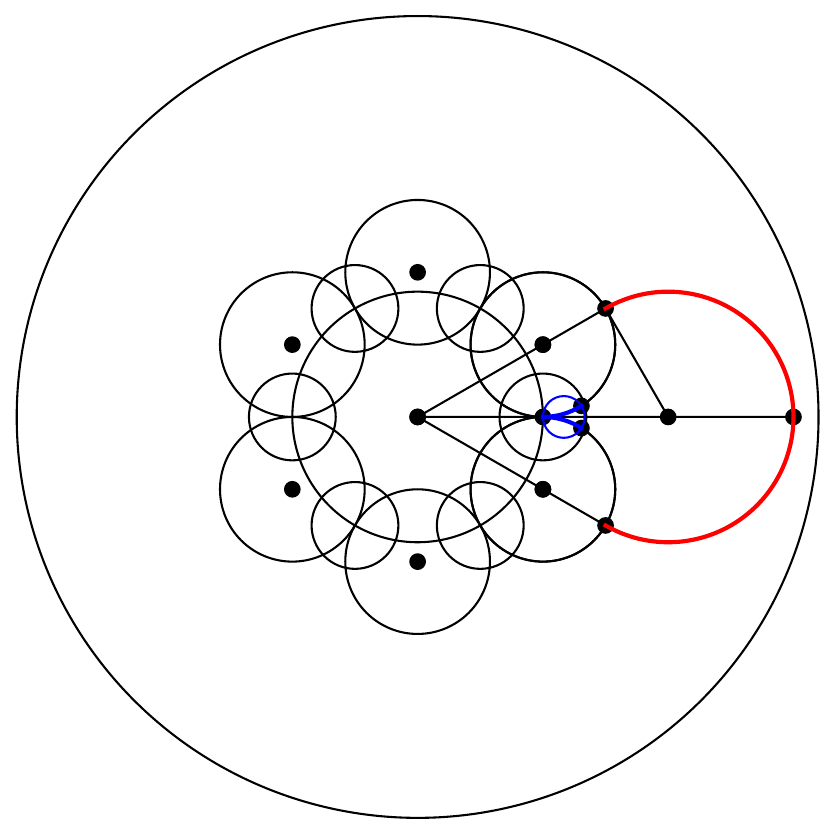}
  }\\
  \subfloat[]{\label{fig:nonvisibleB}
    \includegraphics[width=0.45\textwidth]{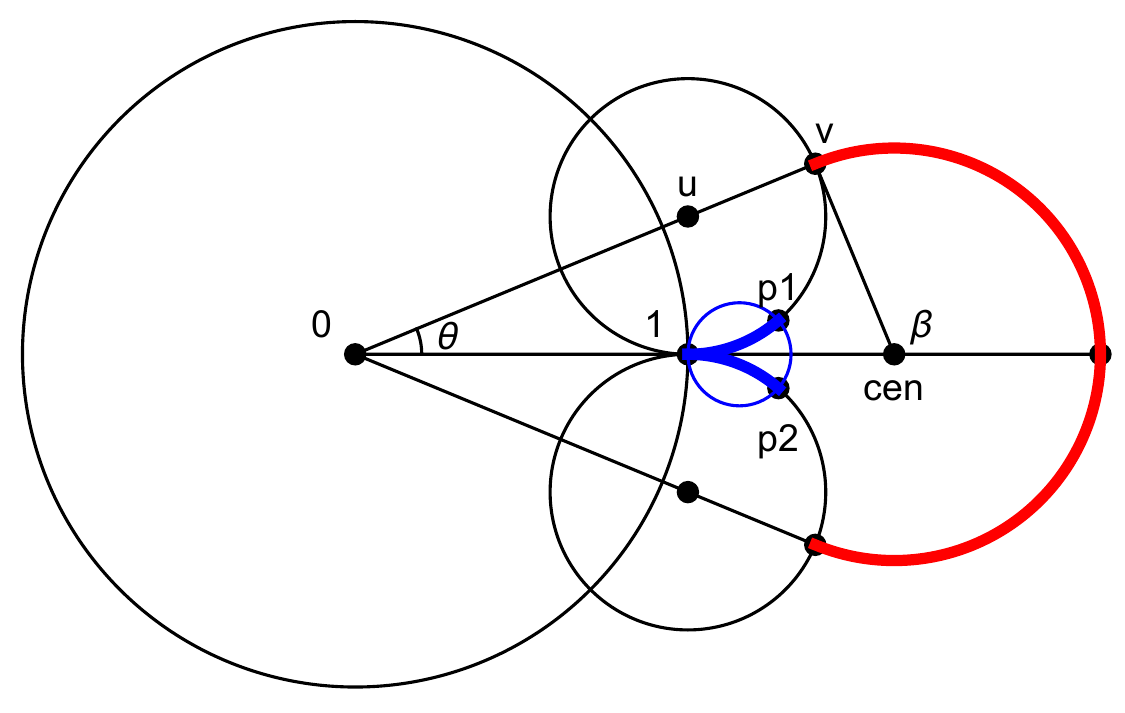}
  }
  \subfloat[]{\label{fig:nonvisibleC}
    \includegraphics[width=0.45\textwidth]{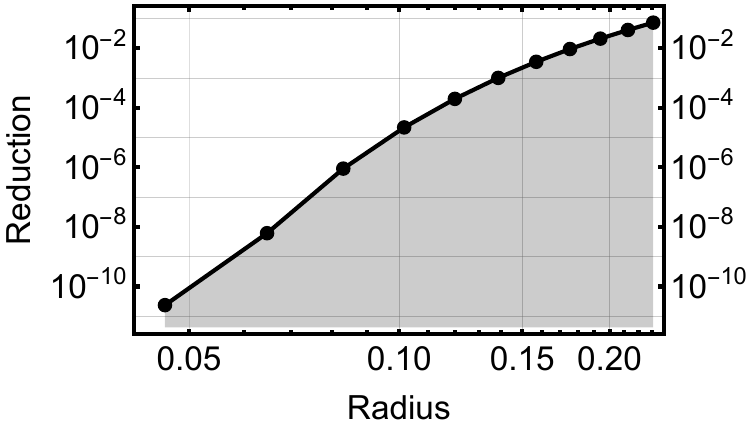}
  }
  \caption{Non-visibility: Top: The full configuration. Bottom left:
  Detailed geometric configuration of the capacity reduction experiment.
  Bottom right: The reduction of the capacity over one segment as function of the radius of the cutting (blue) circle (loglog-plot). Given that the reference capacity $= 2.435704976$, the maximal reduction over the whole set is roughly 3\%.}
  \label{fig:nonvisible}
\end{figure}

It is natural to expect that in the case of Figure~\ref{fig:nonvisibleA} 
certain parts of the compact set $K$ consisting
of the union of mutually tangent circles have very little influence
on the capacity ${ \rm cap}( \mathbb{B}^2,K).$ The parts in question are
the circular arcs close to the contact points of the circles. These parts are
enclosed by the small blue circles. 
Due to symmetry, it is sufficient to consider the effect of area reduction on
one configuration such as the one illustrated in Figure~\ref{fig:nonvisibleB}.
We shall analyse this phenomenon in detail.

\subsection{Conformal mapping of a circular arc triangle onto $\uhp^2$}
We consider a circular arc triangle constructed as follows: Fix $\theta\in(0,\pi/6]$ and let $u=e^{i\theta}/\cos(\theta)$, $v=u+e^{i\theta}{\rm Im}(u)$. Then $\{v\}=(L(0,u)\setminus\B^2)\cap S^1(u,{\rm Im}(u))$. The tangent line to $S^1(u,{\rm Im}(u))$ at the point $v$ intersects the real axis at the point $w={\rm LIS}[0,1,v,v+iv]$ and let $\beta={\rm arg}(v-w)$. Our circular arc triangle consists of three arcs: 
\begin{align*}
S_1=\{z\in S^1(u,{\rm Im}(u))\text{ : }{\rm arg}(z-u)\in[3\pi/2,2\pi+\theta]\},    
\end{align*}
\hspace{3.1cm}$S_2$ equivalent to the complex conjugate of the arc $S_1$, and 
\begin{align*}
S_3=\{z\in S^1(w,|v-w|)\text{ : }{\rm arg}(z-w)\in[2\pi-\beta,2\pi+\beta]\}.    
\end{align*}

Choose a M\"obius transformation $h$ with
\begin{align*}
h(z)=\frac{z+b}{cz-c};\quad h(1)=\infty;\quad h(v)=1;\quad h(\overline{v})=-1.   
\end{align*}
We can solve that
\begin{align*}
b=\frac{v+\overline{v}-2v\overline{v}}{v+\overline{v}-2},\quad
c=\frac{v-\overline{v}}{v+\overline{v}-2}.
\end{align*}
This mapping maps the arcs forming our circular arc triangle as follows:
\begin{align*}
h(S_1)=\{1+it\text{ : }t>0\},\quad
h(S_2)=\{-1+it\text{ : }t>0\},\quad
h(S_3)=(-1,1).
\end{align*}
Consequently, the circular arc triangle is mapped onto the semi-infinite strip
\begin{align*}
T=\{z\text{ : }-1<{\rm Re}(z)<1,\,{\rm Im}(z)>0\}.    
\end{align*}
Finally, the conformal mapping $f(z)=\sin(\pi z/2)$ maps $T$ onto $\uhp^2$ with 
\begin{align*}
f(h(S_1))&=f(\{1+it\text{ : }t>0\})=(\infty,1),\\
f(h(S_2))&=f(\{-1+it\text{ : }t>0\})=(-1,-\infty),\\
f(h(S_3))&=f((-1,1))=(-1,1).
\end{align*}
Thus, to map the circular arc triangle onto $\uhp^2$, the desired mapping is $f\circ h$.

For small values of $s\in(0,|v-w|/3)$, we also consider the sub-arc
\begin{align*}
S_1\cap B^2(ts,s)=\{z\in S_1\text{ : }3\pi/2<{\rm arg}(z-u)<3\pi/2+\theta_2\}    
\end{align*}
of $S_1$, where $\theta_2=2{\rm arctan}(s/{\rm Im}(u))$.


\subsection{Numerical estimation of the loss of capacity}
Using the finite element method one can numerically estimate the loss of capacity due
to reduction of the area. We start by computing ${ \rm cap}( \mathbb{B}^2,K) = 2.435704976.$
Next, we let the radius of the cutting (blue) circles increase and compute the
capacities for every realisation in the sequence. The loss or reduction in capacity
is shown in Figure~\ref{fig:nonvisibleC}. The maximal reduction of  ${ \rm cap}( \mathbb{B}^2,K)$
at cutting radius of $1/4$ is roughly 3\%. As discussed above, these types of configurations
are good candidates for defeaturing.

\section{Conclusions}
Computation of conformal capacity of condensers over complicated domains 
such as the mazes considered here using standard numerical methods 
such as finite element method is expensive. 
The mazes are of importance since they seem have theoretically maximal capacities.
Finding reliable yet easily computable approximations is well motivated in this context. 
Two approximations that are natural are the quasihyperbolic length and perimeter of the compact set.
It has been demonstrated over the set of examples considered that they indeed possess the
right asymptotic properties. Once the necessary analysis is completed, evaluation of the
estimates reduces to evaluating formulae.

Another aspect of the conformal capacity is that in many instances the approximations are amenable
to defeaturing, that is, it is possible to obtain reasonable estimates by removing
geometric features like cusps that complicate for instance meshing. 

\section*{Acknowledgement}
We gratefully acknowledge the memory of Bent Fuglede (1925--2023), who inspired many aspects of this research.
Hakula: This work was supported by the Research Council of Finland (Flagship of Advanced Mathematics for Sensing Imaging and Modelling grant 359181).


\end{document}